\documentclass{amsart}

\usepackage{amssymb}

\allowdisplaybreaks

\newtheorem{theorem}{Theorem}[section]
\newtheorem{lemma}[theorem]{Lemma}

\newtheorem{remark}[theorem]{Remark}

\theoremstyle{definition}

\theoremstyle{definition}

\newcommand{\dem}{\noindent {\bf Proof. }}

\newcommand{\fin}{\hspace*{\fill} $\square$ \vskip0.2cm}

\newcommand{\Word}{\mathcal{W}}

\newcommand{\Free}{\mathbf{F}}

\begin{document}

\title[Multi-indexed $p$-orthogonal sums] {Multi-indexed $p$-orthogonal
sums in non-commutative Lebesgue spaces}

\author[Javier Parcet]
{Javier Parcet}

\address{Universidad Aut\'{o}noma de Madrid and Texas A\&M University}

\email{javier.parcet@uam.es}

\footnote{Partially supported by the Project BFM 2001/0189,
Spain.} \footnote{2000 Mathematics Subject Classification: Primary
46L52. Secondary 05A18.} \footnote{Key words and phrases:
Khintchine inequality, M\"{o}bius inversion, $p$-orthogonal sums.}

\date{}

\begin{abstract}
In this paper we extend a recent Pisier's inequality for
$p$-orthogonal sums in non-commutative Lebesgue spaces. To that
purpose, we generalize the notion of $p$-orthogonality to the
class of multi-indexed families of operators. This kind of
families appear naturally in certain non-commutative Khintchine
type inequalities associated with free groups. Other
$p$-orthogonal families are given by the homogeneous
operator-valued polynomials in the Rademacher variables or the
multi-indexed martingale difference sequences. As in Pisier's
result, our tools are mainly combinatorial.
\end{abstract}

\maketitle

\section*{Introduction}

Let $\mathcal{M}$ be a von Neumann algebra equipped with a
faithful, normal trace $\tau$ satisfying $ \tau(1)=1$ and let us
consider the associated non-commutative Lebesgue space $L_p(\tau)$
for an even integer $p$. Let $\Gamma$ be the product set $\{1,2,
\ldots, n\}^d$ and let $f = (f_{\gamma})_{\gamma \in \Gamma}$ be a
family of operators in $L_p(\tau)$ indexed by $\Gamma$. We shall
say that $f$ is \emph{$p$-orthogonal with $d$ indices} if $$\tau
\big( f_{h(1)}^*f_{h(2)}^{}f_{h(3)}^*f_{h(4)}^{} \cdots
f_{h(p-1)}^*f_{h(p)}^{} \big) = 0$$ whenever the function $h:
\{1,2, \ldots, p\} \rightarrow \Gamma$ has an injective
projection. In other words, whenever the coordinate function
$\pi_k \circ h: \{1,2, \ldots, p\} \rightarrow \{1,2, \ldots, n\}$
is an injective function for some $1 \le k \le d$. Of course, as
it is to be expected, the product above can be replaced by
$$f_{h(1)}^{}f_{h(2)}^* \cdots f_{h(p-1)}^{}f_{h(p)}^*,$$ with no
consequences in the forthcoming results. The case of one index
$d=1$ was already considered by Pisier in \cite{P1}. The main
result in \cite{P1} is the following inequality, which holds for
any $p$-orthogonal family $f_1, f_2, \ldots, f_n$ with one index
$$\Big\| \sum_{k=1}^n f_k \Big\|_{L_p(\tau)} \le \frac{3 \pi}{2} p
\, \max \left\{ \Big\| \Big( \sum_{k=1}^n f_k^* f_k \Big)^{1/2}
\Big\|_{L_p(\tau)}, \Big\| \Big( \sum_{k=1}^n f_k f_k^*
\Big)^{1/2} \Big\|_{L_p(\tau)} \right\}.$$ Some natural examples
of $1$-indexed $p$-orthogonal sequences of operators are the
(non-commutative) \emph{martingale difference sequences}, the
operators associated to a \emph{$p$-dissociate subset} of any
discrete group (via the left regular representation) or a
\emph{free circular family} in Voiculescu's sense \cite{VDN}. In
particular, several relevant inequalities in Harmonic Analysis
such as the Littlewood-Paley inequalities, the (non-commutative)
Burkholder-Gundy inequalities \cite{PX}, or the (non-commutative)
Khintchine inequalities \cite{L,LP} appear as particular cases.
Moreover, it turns out that the combinatorial techniques applied
in \cite{P1} led to the sharp order of growth of the constant
appearing in the non-commutative Burkholder-Gundy inequalities.
For the more general case of $d$ indices, we are interested in
upper bounds for the norm in $L_p(\tau)$ of the sum $$\sum_{\gamma
\in \Gamma} f_{\gamma}.$$ To explain the main result of this
paper, let us introduce some notation. Let $[m]$ be an
abbreviation for the set $\{1,2, \ldots, m\}$. Then, if
$\mathbb{P}_d(2)$ denotes the set of partitions $(\alpha, \beta)$
of $[d]$ into two disjoint subsets (where we allow $\alpha$ and
$\beta$ to be the empty set), we denote by $$\pi_{\alpha}: \Gamma
\rightarrow [n]^{|\alpha|}$$ the canonical projection given by
$\pi_{\alpha}(\gamma) = (i_k)_{k \in \alpha}$ for any $\gamma =
(i_1, \ldots, i_d) \in [n]^d$. Then, if $e_{ij}$ denotes the
natural basis of the Schatten class $S_p$, the sum $$\sum_{\gamma
\in \Gamma} f_{\gamma} \otimes
e_{\pi_{\alpha}(\gamma),\pi_{\beta}(\gamma)}$$ can be understood
as an $L_p(\tau)$-valued matrix with $n^{|\alpha|}$ rows and
$n^{|\beta|}$ columns. In particular, we always obtain an element
of the vector-valued space $L_p(\tau; S_p)$. Our main result can
be stated as follows. Let $p$ be an even integer and let $f =
(f_{\gamma})_{\gamma \in \Gamma}$ be a $p$-orthogonal family in
$L_p(\tau)$ with $d$ indices, then $$\Big\| \sum_{\gamma \in
\Gamma} f_{\gamma} \Big\|_{L_p(\tau)} \le \mathrm{k}_d \,
p^{\frac{d(d+1)}{2}} \max_{(\alpha,\beta) \in \mathbb{P}_d(2)}
\left\{ \Big\| \sum_{\gamma \in \Gamma} f_{\gamma} \otimes
e_{\pi_{\alpha}(\gamma),\pi_{\beta}(\gamma)}
\Big\|_{L_p(\tau;S_p)} \right\}.$$ Here, $\mathrm{k}_d$ denotes an
absolute constant depending only on $d$. Recall that Pisier's
inequality follows from our result for $1$-indexed $p$-orthogonal
sums since $\alpha$ is either $\{1\}$ or the empty set while
$\beta$ is the complement of $\alpha$. The general picture of our
proof follows similar ideas to those in \cite{P1}. Indeed, let
$\Free_n$ be the free group with $n$ generators $g_1, g_2, \ldots,
g_n$ and let $\lambda$ stand for the left regular representation
of $\Free_n$. Then it is easy to check that the family of
operators $$f_{\gamma} = \lambda(g_{i_1}) \otimes \lambda(g_{i_2})
\otimes \cdots \otimes \lambda(g_{i_d}) \qquad \mbox{with} \qquad
\gamma = (i_1,i_2, \ldots, i_d),$$ is $p$-orthogonal with $d$
indices for any even integer $p$. Using the non-commutative
Khintchine inequality for free generators, we show that this
family satisfies the inequality appearing in our main result.
After that, the basic idea is to show that the norm of any
$p$-orthogonal sum with $d$ indices is controlled by the behaviour
of this family. To that aim, we use the same combinatorial
techniques employed in \cite{P1} to obtain a factorization result
which allows us to use H\"{o}lder inequality. Then, the result follows
easily.

In Section \ref{Section-Iteration}, we describe the inequalities
which arise when applying several times the non-commutative
Khintchine inequality for free generators to the family
$\lambda(g_{i_1}) \otimes \cdots \otimes \lambda(g_{i_d})$. These
inequalities will be used in the proof of our result. In Section
\ref{Section-Mobius}, we give a brief summary of results about the
theory of partitions that we shall need in the proof. Section
\ref{Section-Main} is devoted to the proof of the stated
inequality for multi-indexed $p$-orthogonal sums. Section
\ref{Section-Examples} contains two particularly interesting
examples of multi-indexed $p$-orthogonal sums. The first one
generalizes the notion of $p$-dissociate set in a discrete group.
The second one is related to a Burkholder-Gundy type inequality
for multi-indexed martingale difference sequences.

\section{Iterations of the Khintchine inequality}
\label{Section-Iteration}

Let $\Free_n$ be the free group with $n$ generators $g_1,g_2,
\ldots, g_n$. If $\delta_t$ denotes the generic element of the
natural basis of $\ell_2(\Free_n)$, the left regular
representation $\lambda$ of $\Free_n$ is defined by the relation
$$\lambda(t_1) \delta_{t_2} = \delta_{t_1t_2}.$$ The reduced
$C^*$-algebra $C_{\lambda}^*(\Free_n)$ is defined as the
$C^*$-algebra generated in $\mathcal{B}(\ell_2(\Free_n))$ by the
operators $\lambda(t)$ when $t$ runs over $\Free_n$. Let us denote
by $\tau$ the standard trace on $C_{\lambda}^*(\Free_n)$ defined
by $\tau(x) = \langle x \delta_e, \delta_e \rangle$, where $e$
denotes the identity element of $\Free_n$. Then, we construct the
non-commutative Lebesgue space $L_p(\tau)$ in the usual way and
consider the subspace $\Word_p(n)$ of $L_p(\tau)$ generated by the
operators $\lambda(g_1), \lambda(g_2), \ldots, \lambda(g_n)$. The
next result was proved by Haagerup and Pisier in \cite{HP} when
$p=\infty$ and extended to any exponent $2 \le p \le \infty$ in
\cite{P2}.

\begin{lemma} \label{Lemma-Iteration}
Let $a_1, a_2, \ldots, a_n$ be a family of operators in some
non-commutative Lebesgue space $L_p(\varphi)$. The following
equivalence of norms holds for $2 \le p \le \infty$, $$\Big\|
\sum_{k=1}^n a_k \otimes \lambda(g_k) \Big\|_{L_p(\varphi \otimes
\tau)} \simeq \max \left\{ \Big\| \sum_{k=1}^n a_k \otimes e_{1k}
\Big\|_{L_p(\varphi; R_p^n)} , \Big\| \sum_{k=1}^n a_k \otimes
e_{k1} \Big\|_{L_p(\varphi; C_p^n)} \right\} \! .$$ In fact, the
linear map $u: R_p^n \cap C_p^n \rightarrow \Word_p(n)$ defined by
$$u(e_{1k} \oplus e_{k1}) = \lambda(g_k),$$ is a complete
isomorphism with $\|u\|_{cb} \le 2$ and completely contractive
inverse.
\end{lemma}

The row and column Hilbert spaces $R_p^n$ and $C_p^n$ are defined
as the operator spaces generated by $\{e_{1j} \, | \ 1 \le j \le
n\}$ and $\{e_{i1} \, | \ 1 \le i \le n\}$ respectively in $S_p$.
Now, let us consider the group product $\mathrm{G}_d = \Free_n
\times \Free_n \times \cdots \times \Free_n$ with $d$ factors. The
left regular representation $\lambda_d$ of $\mathrm{G}_d$ has the
form $$\lambda_d(t_1,t_2, \ldots, t_d) = \lambda(t_1) \otimes
\lambda(t_2) \otimes \cdots \otimes \lambda(t_d).$$ Hence, the
reduced $C^*$-algebra $C_{\lambda_d}^*(\mathrm{G}_d)$ is endowed
with the trace $\tau_d = \tau \otimes \tau \otimes \cdots \otimes
\tau$ with $d$ factors. This allows us to consider the
non-commutative space $L_p(\tau_d)$ for any $1 \le p \le \infty$.
Then we define the space $\Word_p(n)^{\otimes d}$ to be the
subspace of $L_p(\tau_d)$ generated by the family of operators
$$\lambda(g_{i_1}) \otimes \lambda(g_{i_2}) \otimes \cdots \otimes
\lambda(g_{i_d}).$$

The aim of this section is to describe the operator space
structure of $\Word_p(n)^{\otimes d}$ as a subspace of
$L_p(\tau_d)$ for the exponents $2 \le p \le \infty$. This
operator space structure has been already described in
\cite[Section 9.8]{P2}, but here we shall give a more detailed
exposition. As it was pointed out in \cite{P2}, the case $1 \le p
\le 2$ follows easily by duality. However, we shall not write the
explicit inequalities in that case since we are not using them and
the notation is considerably more complicated. If we apply
repeatedly Lemma \ref{Lemma-Iteration} to the sum
$$\mathcal{S}_d(a) = \sum_{i_1, \ldots, i_d=1}^n a_{i_1 i_2 \cdots
i_d} \otimes \lambda(g_{i_1}) \otimes \lambda(g_{i_2}) \otimes
\cdots \otimes \lambda(g_{i_d}) \in L_p(\varphi \otimes \tau_d),$$
then we easily get $$\|\mathcal{S}_d(a)\|_{L_p(\varphi \otimes
\tau_d)} \le 2^d \max \left\{ \Big\| \sum_{i_1, \ldots, i_d =1}^n
a_{i_1 \cdots i_d} \otimes \xi_1(i_1) \otimes \cdots \otimes
\xi_d(i_d) \Big\|_{L_p(\varphi; S_p)} \ \right\},$$ where the
maximum runs over all possible ways to choose the functions
$\xi_1, \xi_2, \ldots, \xi_d$ among $\xi_k(\cdot) = e_{\cdot 1}$
and $\xi_k(\cdot) = e_{1 \cdot}$. That is, each function $\xi_k$
can take values either in the space $R_p^n$ or in the space
$C_p^n$. For a given selection of $\xi_1, \xi_2, \ldots, \xi_d$ we
split up these functions into two sets, one made up of the
functions taking values in $R_p^n$ and the other taking values in
$C_p^n$. More concretely, let us consider the sets
$$\begin{array}{l} \mathrm{R}_{\xi} = \big\{ k \, | \ \xi_k(i) =
e_{1i} \big\}, \\ \mathrm{C}_{\xi} = \big\{ k \, | \ \xi_k(i) =
e_{i1} \big\}. \end{array}$$ Then, if $\mathrm{C}_{\xi}$ has $s$
elements, the sum $$\sum_{i_1, \ldots, i_d =1}^n a_{i_1 \cdots
i_d} \otimes \xi_1(i_1) \otimes \cdots \otimes \xi_d(i_d)$$ can be
regarded as a $n^s \times n^{d-s}$ matrix with entries in
$L_p(\varphi)$. Now, using the notation already presented in the
Introduction, we express the inequality above in a much more
understandable way. Namely, we have
\begin{equation} \label{Equation-Iteration}
\|\mathcal{S}_d(a)\|_{L_p(\varphi \otimes \tau_d)} \le 2^d
\max_{(\alpha,\beta) \in \mathbb{P}_d(2)} \left\{ \Big\|
\sum_{\gamma \in \Gamma} a_{\gamma}^{} \otimes
e_{\pi_{\alpha}(\gamma),\pi_{\beta}(\gamma)} \Big\|_{L_p(\varphi;
S_p)} \ \right\}.
\end{equation}

\begin{remark} \label{Remark-Converse}
\textnormal{By the same arguments, the converse of
(\ref{Equation-Iteration}) holds with constant 1.}
\end{remark}

\section{M\"{o}bius inversion}
\label{Section-Mobius}

Given a positive integer $m$, we denote by $\mathbb{P}_m$ the
lattice of partitions of the set $[m] = \{1,2, \ldots, m\}$. If
$\rho$ and $\sigma$ are elements of $\mathbb{P}_m$, we shall write
$\rho \le \sigma$ when every block of $\rho$ is contained in some
block of $\sigma$. The minimal and maximal elements of
$\mathbb{P}_m$ with respect to this partial order are denoted by
$\dot{0}$ and $\dot{1}$ respectively. That is, $\dot{0}$ stands
for the partition into $m$ singletons and $\dot{1}$ coincides with
$\{[m]\}$. The M\"{o}bius function $\mu$ is a complex-valued function
defined on the set of pairs of partitions $(\rho, \sigma)$ in
$\mathbb{P}_m \times \mathbb{P}_m$ satisfying $\rho \le \sigma$.
The following Lemma summarizes the main properties of this
function that we shall use below.

\begin{lemma} \label{Lemma-Mobius}
Let us consider a pair of functions $\Phi: \mathbb{P}_m
\rightarrow V$ and $\Psi: \mathbb{P}_m \rightarrow V$ taking
values in some vector space $V$. Then the following implication
holds $$\Psi(\rho) = \sum_{\sigma \ge \rho} \Phi(\sigma) \
\Rightarrow \ \Phi(\rho) = \sum_{\sigma \ge \rho} \mu(\rho,
\sigma) \Psi(\sigma).$$ Besides, the M\"{o}bius function satisfies the
following identities
\begin{itemize}
\item $\displaystyle \sum_{\sigma \in \mathbb{P}_m}
|\mu(\dot{0},\sigma)| = m!.$
\item For any $\sigma > \dot{0}$, we have $\displaystyle
\sum_{\dot{0} \le \rho \le \sigma} \mu(\rho, \sigma) = 0$.
\end{itemize}
\end{lemma}

For a more detailed exposition of these topics we refer the reader
to \cite{A}. Now, let $p$ be an even integer and let $\varphi: E_1
\times \cdots \times E_p \rightarrow V$ be a multilinear map
defined on certain vector spaces $E_1, E_2, \ldots, E_p$ and
taking values in the vector space $V$. For each $1 \le s \le p$ we
consider elements $f_{\gamma}(s) \in E_s$ indexed by $\Gamma$.
Then, we define the sums $$\mathrm{F}_s = \sum_{\gamma \in \Gamma}
f_{\gamma}(s) \in E_s.$$ Clearly we have $$\varphi(\mathrm{F}_1,
\mathrm{F}_2 ,\ldots, \mathrm{F}_p) = \sum_h \varphi
\big(f_{h(1)}(1), f_{h(2)}(2), \ldots, f_{h(p)}(p) \big),$$ where
the sum runs over the set of functions $h: \{1,2, \ldots, p \}
\rightarrow \Gamma$. Now, for any such function $h$ and for each
$1 \le k \le d$, we consider the partition $\sigma_k(h) \in
\mathbb{P}_p$ associated to the coordinate function $\pi_k \circ
h$. In other words, given $1 \le r,s \le p$ we have the following
characterization $$r \sim s \, (\mbox{mod} \, \sigma_k(h))
\Leftrightarrow \pi_k(h(r)) = \pi_k(h(s)),$$ where $\sim \,
(\mbox{mod} \, \sigma)$ means belonging to the same block of
$\sigma$. Let us also consider the $d$-tuple $\delta(h) =
(\sigma_1(h), \sigma_2(h), \ldots, \sigma_d(h))$ in the product
$\mathbf{P}(p,d) = \mathbb{P}_p \times \cdots \times \mathbb{P}_p$
with $d$ factors. Then we can write $$\varphi(\mathrm{F}_1,
\ldots, \mathrm{F}_p) = \sum_{\eta \in \mathbf{P}(p,d)}
\Phi(\eta),$$ where $\Phi: \mathbf{P}(p,d) \rightarrow V$ has the
form $$\Phi(\eta) = \sum_{h: \, \delta(h) = \eta} \varphi \big(
f_{h(1)}(1), \ldots, f_{h(p)}(p) \big).$$ Now, if $\eta = (\rho_1,
\ldots , \rho_d)$ we shall write $\eta \sim \dot{0}$ whenever
$\rho_k = \dot{0}$ for some $1 \le k \le d$. Then, we obtain the
following decomposition
\begin{equation} \label{Equation-Decomposition}
\varphi(\mathrm{F}_1, \ldots, \mathrm{F}_p) = \sum_{\eta \sim
\dot{0}} \Phi(\eta) + \sum_{\rho_1 > \dot{0}} \cdots \sum_{\rho_d
> \dot{0}} \Phi(\eta).
\end{equation}
Similarly, the expression $h \sim \dot{0}$ will denote the
existence of some $1 \le k \le d$ such that $\sigma_k(h) =
\dot{0}$. In other words, $h \sim \dot{0}$ whenever $h$ has an
injective projection. Then, since $\delta(h) = (\sigma_1(h),
\ldots, \sigma_d(h))$, we have
\begin{equation} \label{Equation-Injective}
\sum_{\eta \sim \dot{0}} \Phi(\eta) = \sum_{h \sim \dot{0}}
\varphi \big( f_{h(1)}(1), \ldots, f_{h(p)}(p) \big).
\end{equation}
For the second sum in (\ref{Equation-Decomposition}), we define
$$\Psi_d(\rho_1, \rho_2, \ldots, \rho_{d}) = \sum_{\sigma_d \ge
\rho_d} \Phi(\rho_1, \rho_2, \ldots, \rho_{d-1} | \sigma_d).$$
Then, if we fix $\rho_1, \rho_2, \ldots, \rho_{d-1}$, we can apply
Lemma \ref{Lemma-Mobius} to obtain
\begin{eqnarray*}
\sum_{\rho_d > \dot{0}} \Phi(\rho_1, \rho_2, \ldots, \rho_d) & = &
\sum_{\rho_d > \dot{0}} \Big( \sum_{\sigma_d \ge \rho_d}
\mu(\rho_d,\sigma_d) \Psi_d(\rho_1, \ldots, \rho_{d-1} | \sigma_d)
\Big) \\ & = & \sum_{\sigma_d > \dot{0}} \Psi_d(\rho_1, \ldots,
\rho_{d-1} | \sigma_d) \sum_{\dot{0} < \rho_d \le \sigma_d}
\mu(\rho_d,\sigma_d) \\ & = & \sum_{\sigma_d
> \dot{0}} (- \mu(\dot{0},\sigma_d)) \Psi_d(\rho_1, \ldots, \rho_{d-1}
| \sigma_d)
\end{eqnarray*}
Similarly, we define
\begin{eqnarray*}
\Psi_{d-1}(\rho_1, \rho_2, \ldots, \rho_{d-1} | \sigma_d) & = &
\sum_{\sigma_{d-1} \ge \rho_{d-1}} \Psi_d(\rho_1, \rho_2, \ldots,
\rho_{d-2} | \sigma_{d-1}, \sigma_d), \\ \Psi_{d-2}(\rho_1,
\rho_2, \ldots | \sigma_{d-1}, \sigma_d) & = & \sum_{\sigma_{d-2}
\ge \rho_{d-2}} \Psi_{d-1} (\rho_1, \rho_2, \ldots | \sigma_{d-2},
\sigma_{d-1}, \sigma_d),
\end{eqnarray*}
and so on until $$\Psi_{1}(\rho_1 | \sigma_2, \ldots, \sigma_d) =
\sum_{\sigma_1 \ge \rho_1} \Psi_2(\sigma_1, \sigma_2, \ldots,
\sigma_d).  \ \qquad \qquad \null$$ Then, applying Lemma
\ref{Lemma-Mobius} as above, we have for $1 \le k \le d-1$
$$\sum_{\rho_k > \dot{0}}^{\null} \Psi_{k+1}(\rho_1, \ldots |
\sigma_{k+1}, \ldots, \sigma_d) = - \sum_{\sigma_k > \dot{0}}
\mu(\dot{0}, \sigma_k) \Psi_k(\rho_1, \ldots | \sigma_k, \ldots,
\sigma_d).$$ Putting all together, we get
\begin{eqnarray} \label{Equation-Non-Injective}
\sum_{\rho_1 > \dot{0}} \cdots \sum_{\rho_d > \dot{0}} \Phi(\eta)
& = & (-1)^d \sum_{\sigma_1 > \dot{0}} \cdots \sum_{\sigma_d >
\dot{0}} \Big[ \prod_{k=1}^d \mu(\dot{0},\sigma_k) \Big]
\Psi_1(\sigma_1, \ldots, \sigma_d),
\end{eqnarray}
where the function $\Psi_1$ can be easily rewritten as
\begin{equation} \label{Equation-Psi1}
\Psi_1(\sigma_1, \ldots, \sigma_d) = \sum_{\begin{array}{c}
\mbox{\footnotesize{$h: \, \sigma_k(h) \ge \sigma_k$}} \\
\mbox{\footnotesize{$1 \le k \le d$}}
\end{array}}^{\null} \varphi \big( f_{h(1)}(1), f_{h(2)}(2),
\ldots, f_{h(p)}(p) \big).
\end{equation}
In summary, looking at (\ref{Equation-Decomposition}),
(\ref{Equation-Injective}), (\ref{Equation-Non-Injective}) and
(\ref{Equation-Psi1}) we have the following result.

\begin{lemma} \label{Lemma-Injective}
The following identity holds
\begin{eqnarray*}
\varphi(\mathrm{F}_1, \ldots, \mathrm{F}_p) & = & \sum_{h \sim
\dot{0}} \varphi \big( f_{h(1)}(1), \ldots, f_{h(p)}(p) \big)
\\ & + & (-1)^d \sum_{\sigma_1 > \dot{0}} \cdots \sum_{\sigma_d >
\dot{0}} \Big[ \prod_{k=1}^d \mu(\dot{0},\sigma_k) \Big]
\Psi(\sigma_1, \ldots, \sigma_d),
\end{eqnarray*} where $\Psi$ has the following form
$$\Psi(\sigma_1, \ldots, \sigma_d) = \sum_{\begin{array}{c}
\mbox{\footnotesize{$h: \, \sigma_k(h) \ge \sigma_k$}} \\
\mbox{\footnotesize{$1 \le k \le d$}}
\end{array}} \varphi \big( f_{h(1)}(1), f_{h(2)}(2), \ldots,
f_{h(p)}(p) \big).$$
\end{lemma}

\section{Proof of the main result}
\label{Section-Main}

In this section we shall prove the result stated below. We start
by factorizing the sum which defines the function $\Psi$ above.
This will allow us to show that the behaviour of any
$p$-orthogonal sum with $d$ indices is majorized by the estimates
obtained in Section \ref{Section-Iteration}, with the aid of
non-commutative Khintchine inequalities.

\begin{theorem} \label{Theorem-p-Orthogonality}
If $f = (f_{\gamma})_{\gamma \in \Gamma}$ is $p$-orthogonal in
$L_p(\tau)$ with $d$ indices, then $$\Big\| \sum_{\gamma \in
\Gamma} f_{\gamma} \Big\|_{L_p(\tau)} \le \mathrm{k}_d \,
p^{\frac{d(d+1)}{2}} \max_{(\alpha,\beta) \in \mathbb{P}_d(2)}
\left\{ \Big\| \sum_{\gamma \in \Gamma} f_{\gamma} \otimes
e_{\pi_{\alpha}(\gamma),\pi_{\beta}(\gamma)}
\Big\|_{L_p(\tau;S_p)} \right\}.$$
\end{theorem}

\subsection{Factorization of $\Psi$}

Let $\mathcal{M}$ be a von Neumann algebra equipped with a
faithful normal trace $\tau$ satisfying $\tau(1)=1$ and let $p$ be
an even integer. Following the notation above, we shall take in
what follows $E_s = L_p(\tau)$ for all $1 \le s \le p$ and the
multilinear map $\varphi$ will be replaced by the trace $\tau$
acting on a product of $p$ operators in $L_p(\tau)$. That is,
$f=(f_{\gamma})_{\gamma \in \Gamma}$ is assumed to be
$p$-orthogonal in $L_p(\tau)$ with $d$ indices and we have
$$\varphi \big( f_{h(1)}(1), f_{h(2)}(2), \ldots, f_{h(p)}(p)
\big) = \tau \big( f_{h(1)}(1) f_{h(2)}(2) \cdots f_{h(p)}(p)
\big),$$ where $$f_{h(s)}(s) = \left\{
\begin{array}{ll} f_{h(s)}^* & \mbox{if $s$ is odd,} \\ f_{h(s)} &
\mbox{if $s$ is even.}
\end{array} \right.$$ The aim now is to factorize the sum
$$\Psi(\sigma_1, \ldots, \sigma_d) = \sum_{\begin{array}{c}
\mbox{\footnotesize{$h: \, \sigma_k(h) \ge \sigma_k$}} \\
\mbox{\footnotesize{$1 \le k \le d$}} \end{array}}^{\null} \tau
\big( f_{h(1)}(1) f_{h(2)}(2) \cdots f_{h(p)}(p) \big).$$ We shall
need below the following version of Fell's absorption principle.

\vskip5pt

\noindent \textbf{Absorption Principle in $L_p$.} \emph{Given a
discrete group $\mathrm{G}$, let us denote by
$\lambda_{\mathrm{G}}$ the left regular representation of
$\mathrm{G}$ and by $\tau_{\mathrm{G}}$ the associated trace on
the reduced $C^*$-algebra of $\mathrm{G}$. Then, given any other
unitary representation $\pi: \mathrm{G} \rightarrow
\pi(\mathrm{G})''$, the following representations are unitarily
equivalent $$\lambda_{\mathrm{G}} \otimes \pi \simeq
\lambda_{\mathrm{G}} \otimes 1,$$ where $1$ stands for the trivial
representation of $\mathrm{G}$ in $\pi(\mathrm{G})''$. Let us
consider any faithful normalized trace $\psi$ on
$\pi(\mathrm{G})''$. Then, given any finitely supported function
$a: \mathrm{G} \rightarrow L_p(\varphi)$, the following equality
holds for $1 \le p \le \infty$} $$\Big\| \sum_{t \in \mathrm{G}}
a(t) \otimes \lambda_{\mathrm{G}}(t) \otimes \pi(t)
\Big\|_{L_p(\varphi \otimes \tau_{\mathrm{G}} \otimes \psi)} =
\Big\| \sum_{t \in \mathrm{G}} a(t) \otimes
\lambda_{\mathrm{G}}(t) \Big\|_{L_p(\varphi \otimes
\tau_{\mathrm{G}})}.$$

\vskip5pt

\dem See Proposition 8.1 of \cite{P2} for the first part and
\cite{PP} for the second. \fin

\begin{lemma} \label{Lemma-Factorization}
Let $\sigma_1, \sigma_2, \ldots, \sigma_d$ be a family of
partitions in $\mathbb{P}_p$ different from $\dot{0}$. If we are
given $0 \le q \le d$, let $\mathrm{B}_q$ be the set of elements
$s$ in $\{1,2, \ldots, p\}$ being a singleton exactly in $q$
partitions among $\sigma_1, \sigma_2, \ldots, \sigma_d$. Then,
there exists a discrete group $\mathrm{G}$ and a family
$\mathsf{F}_1, \mathsf{F}_2, \ldots, \mathsf{F}_p$ in
$L_p(\tau_{\mathrm{G}} \otimes \tau)$ satisfying
$$\|\mathsf{F}_s\|_p \le \mathrm{k}_d \, p^{\frac{q(q+1)}{2}}
\Big\| \sum_{i_1, \ldots, i_d = 1}^n \lambda(g_{i_1}) \otimes
\cdots \otimes \lambda(g_{i_d}) \otimes f_{i_1 \cdots i_d}
\Big\|_{L_p(\tau_d \otimes \tau)}$$ for each $s \in \mathrm{B}_q$
whenever $0 \le q < d$ and also $$\|\mathsf{F}_s\|_p = \Big\|
\sum_{\gamma \in \Gamma} f_{\gamma} \Big\|_{L_p(\tau)}$$ for each
$s \in \mathrm{B}_d$. Moreover, we have
\begin{equation} \label{Equation-Factorization}
\sum_{\begin{array}{c} \mbox{\footnotesize{$h: \, \sigma_k(h) \ge
\sigma_k$}} \\ \mbox{\footnotesize{$1 \le k \le d$}}
\end{array}}^{\null} \tau \big( f_{h(1)}(1) f_{h(2)}(2) \cdots
f_{h(p)}(p) \big) = (\tau_{\mathrm{G}} \otimes \tau) (\mathsf{F}_1
\mathsf{F}_2 \cdots \mathsf{F}_p).
\end{equation}
\end{lemma}

\begin{remark} \label{Remark-Induction}
\textnormal{As we have pointed out, Theorem
\ref{Theorem-p-Orthogonality} was already proved in \cite{P1} for
$1$-indexed families. In particular, we can assume that Theorem
\ref{Theorem-p-Orthogonality} holds for any $k$-indexed family
whenever $1 \le k \le d-1$ and prove Theorem
\ref{Theorem-p-Orthogonality} by induction. In the proof of Lemma
\ref{Lemma-Factorization}, we shall need to use this induction
hypothesis.}
\end{remark}

\begin{remark}
\textnormal{From now on, $\mathrm{k}_d$ might change from one
instance to another.}
\end{remark}

\dem Let us consider an integer $2 \le m \le p$. As it is
customary, we write $\tau_{m-1}$ for the standard trace associated
to the reduced $C^*$-algebra of the group product $\Free_n \times
\Free_n \times \cdots \times \Free_n$ with $m-1$ factors. Then,
for each $1 \le i \le n$, we consider the following family
$\xi_1(i), \xi_2(i), \ldots, \xi_m(i)$ of operators in
$L_p(\tau_{m-1})$
\begin{eqnarray*}
\xi_1(i) & = & \lambda(g_i)^* \otimes 1 \otimes 1 \otimes \cdots
\otimes 1 \otimes 1, \\ \xi_2(i) & = & \lambda(g_i) \otimes
\lambda(g_i)^* \otimes 1 \otimes \cdots \otimes 1, \\\xi_3(i) & =
& 1 \otimes \lambda(g_i) \otimes \lambda(g_i)^* \otimes \cdots
\otimes 1, \\ & \cdots & \\ \xi_{m-1}(i) & = & 1 \otimes \cdots
\otimes 1 \otimes \lambda(g_i) \otimes \lambda(g_i)^*, \\ \xi_m(i)
& = & 1 \otimes \cdots \otimes 1 \otimes \lambda(g_i).
\end{eqnarray*}
Given $g: \{1,2, \ldots, m\} \rightarrow \{1,2, \ldots,n\}$, this
family has the following property
\begin{equation} \label{Equation-Non-Constant}
\tau_{m-1} \big( \xi_1(g(1)) \cdots \xi_m(g(m)) \big) = \left\{
\begin{array}{ll} 1 & \mbox{if $g$ is constant,} \\ 0 & \mbox{if
$g$ is non-constant}.
\end{array} \right.
\end{equation}
Let us make explicit the blocks of the partitions $\sigma_1,
\sigma_2, \ldots, \sigma_d$ by $$\sigma_k = \Big\{
\mathrm{A}_{kj_k} \, \Big| \ 1 \le j_k \le \mathrm{m}_k \Big\}.$$
Now we fix $\sigma_k$ and, for each $\mathrm{A}_{kj_k}$ with
cardinality $m_{j_k} > 1$, we construct the family $\Pi(i,j_k) =
\{ \xi_1(i,j_k), \xi_2(i,j_k), \ldots, \xi_{m_{j_k}}(i,j_k) \}$ in
$L_p(\tau_{m_{j_k}-1})$ as above. Notice that $1 \le i \le n$ and
$1 \le j_k \le \mathrm{m}_k$. If the set $\mathrm{A}_{kj_k}$ has
only one element, we take $\Pi(i,j_k) = \{ \xi_1(i,j_k) \}$ with
$\xi_1(i,j_k) = 1$. Then we consider the following families of
$\mathrm{m}_k$-fold tensor products
\begin{eqnarray*}
\Sigma(i,1) & = & \Pi(i,1) \otimes 1 \otimes 1 \otimes \cdots
\otimes 1, \\ \Sigma(i,2) & = & 1 \otimes \Pi(i,2) \otimes 1
\otimes \cdots \otimes 1, \\ & \cdots & \\ \Sigma(i,\mathrm{m}_k)
& = & 1 \otimes 1 \otimes \cdots \otimes 1 \otimes
\Pi(i,\mathrm{m}_k).
\end{eqnarray*}
Here, the $r$-th \lq$1$\rq${}$ in $\Sigma(i,j_k)$ denotes the
identity operator in $L_p(\tau_{m_r - 1})$. Recall that, fixed $1
\le i \le n$, each $\Sigma(i, j_k)$ is an \emph{ordered family}
with $m_{j_k}$ elements. On the other hand, for each $1 \le s \le
p$, there exist a unique set of indices $j_1(s), j_2(s), \ldots,
j_d(s)$ such that $s$ belongs to the corresponding blocks of
$\sigma_1, \sigma_2, \ldots, \sigma_d$. In other words, we pick up
the indices $j_k(s)$ satisfying $$s \in \bigcap_{k=1}^d
\mathrm{A}_{kj_k(s)}.$$ This allows us to consider the family of
operators $$\Lambda(\gamma,s) = \bigotimes_{k=1}^d
\Sigma(i_k,j_k(s)),$$ where $i_k = \pi_k(\gamma)$ and $\gamma \in
\Gamma$. Now we select an element of $\Lambda(\gamma,s)$ as
follows. If $s$ is the $r_1$-th element in the block
$\mathrm{A}_{1j_1(s)}$, then we pick up the $r_1$-th operator in
the family $\Sigma(i_1,j_1(s))$. Let us denote it by
$x_{1s}(i_1)$. Similarly, if $s$ is the $r_2$-th element in
$\mathrm{A}_{2j_2(s)}$, we pick up the $r_2$-th operator in
$\Sigma(i_2,j_2(s))$, say $x_{2s}(i_2)$. We iterate this process
to get an element $$x_{1s}(i_1) \otimes x_{2s}(i_2) \otimes \cdots
\otimes x_{ds}(i_d) \in \Lambda(\gamma,s).$$ Then we define,
$$\mathsf{F}_s = \sum_{\gamma \in \Gamma} \Big( \bigotimes_{k=1}^d
x_{ks}(\pi_k(\gamma)) \Big) \otimes f_{\gamma}(s).$$ Clearly,
there exists a collection of discrete groups $\mathrm{G}_1,
\mathrm{G}_2, \ldots, \mathrm{G}_d$ (all of them being direct
products of $\Free_n$) such that $\mathsf{F}_s \in
L_p(\tau_{\mathrm{G}} \otimes \tau)$ with $\mathrm{G} =
\mathrm{G}_1 \times \cdots \times \mathrm{G}_d$. Let us check that
identity (\ref{Equation-Factorization}) holds. Notice that
$$(\tau_{\mathrm{G}} \otimes \tau) (\mathsf{F}_1 \mathsf{F}_2
\cdots \mathsf{F}_p) = \sum_{\gamma_1, \ldots, \gamma_p \in
\Gamma} \prod_{k=1}^d \tau_{\mathrm{G}_k} \Big( \prod_{s=1}^p
x_{ks}(\pi_k(\gamma_s)) \Big) \tau(f_{\gamma_1}(1) \cdots
f_{\gamma_p}(p)).$$ Recalling the definition of $x_{ks}$ and
property (\ref{Equation-Non-Constant}), it can be checked that
$$\prod_{k=1}^d \tau_{\mathrm{G}_k} \Big( \prod_{s=1}^p
x_{ks}(\pi_k(\gamma_s)) \Big)$$ is $1$ when the condition $j_k(s)
= j_k(s') \Rightarrow \pi_k(\gamma_s) = \pi_k(\gamma_{s'})$ holds
for $k=1,2, \ldots, d$ and is $0$ otherwise. In particular,
identity (\ref{Equation-Factorization}) follows. Now we look at
the norm of $\mathsf{F}_s$ in $L_p(\tau_{\mathrm{G}} \otimes
\tau)$. First assume that $s \in \mathrm{B}_d$. That is, $s$ is a
singleton of $\sigma_k$ for every $k = 1,2, \ldots, d$. Then
$$\mathsf{F}_s = \sum_{\gamma \in \Gamma} 1 \otimes f_{\gamma}(s)
= \Big( \sum_{\gamma \in \Gamma} 1 \otimes f_{\gamma}
\Big)^{(*)}$$ where $(*)$ is $*$ when $s$ is odd and $1$
otherwise. Therefore the stated assertion follows. Finally, assume
that $s \in \mathrm{B}_q$ with $q < d$. If $q = 0$ our estimation
for the norm of $\mathsf{F}_s$ is easy. Namely, a quick inspection
of the definition of $\mathsf{F}_s$ allows us to write
$$\mathsf{F}_s \simeq \sum_{i_1, \ldots, i_d =1}^n \Big(
\bigotimes_{k=1}^d \chi(i_k) \Big) \otimes f_{i_1 \cdots
i_d}(s),$$ where $\chi(i_k)$ can be either $\lambda(g_{i_k})^*$ or
$\lambda(g_{i_k}) \otimes \lambda(g_{i_k})^*$ or
$\lambda(g_{i_k})$. However, by Fell's absorption principle these
terms are unitarily equivalent. In other words, in this particular
case we obtain an equality $$\|\mathsf{F}_s\|_p = \Big\|
\sum_{i_1, \ldots, i_d =1}^n \lambda(g_{i_1}) \otimes \cdots
\otimes \lambda(g_{i_d}) \otimes f_{i_1 \cdots i_d}(s)
\Big\|_{L_p(\tau_d \otimes \tau)}.$$ Notice that the dependence on
$s$ on the right can be ignored since the two possible expressions
that come out (for $s$ odd and $s$ even) turn out to be equal. It
remains to check the cases $0 < q < d$. For simplicity of
notation, we assume that $s$ is a singleton in $\sigma_1,
\sigma_2, \ldots, \sigma_q$. As we shall see, the general case can
be proved in a similar way. Then, again by Fell's absorption
principle, we have $$\|\mathsf{F}_s\|_p = \Big\| \sum_{i_1,
\ldots, i_d =1}^n \lambda(g_{i_{q+1}}) \otimes \cdots \otimes
\lambda(g_{i_d}) \otimes f_{i_1 \cdots i_d}(s)
\Big\|_{L_p(\tau_{d-q} \otimes \tau)}.$$ Applying the iteration of
Khintchine inequality described in (\ref{Equation-Iteration}), we
have $$\|\mathsf{F}_s\|_p \le \mathrm{k}_d \, \max_{(\alpha,
\beta) \in \mathbb{P}_{d-q}(2)} \left\{ \Big\| \sum_{\nu \in
[n]^{d-q}} \Big[ \sum_{\zeta \in [n]^q} f_{\zeta, \nu}(s) \Big]
\otimes e_{\pi_{\alpha} (\nu),\pi_{\beta}(\nu)} \Big\|_{L_p(\tau;
S_p)} \right\}.$$ The sum on the right can be rewritten as follows
\begin{eqnarray*}
\sum_{\nu \in [n]^{d-q}} \Big[ \sum_{\zeta \in [n]^q} f_{\zeta,
\nu}(s) \Big] \otimes e_{\pi_{\alpha} (\nu),\pi_{\beta}(\nu)} & =
& \sum_{\zeta \in [n]^q} \Big[ \sum_{\nu \in [n]^{d-q}} f_{\zeta,
\nu}(s) \otimes e_{\pi_{\alpha}(\nu),\pi_{\beta}(\nu)} \Big]
\\ & = & \sum_{\zeta \in [n]^q} f_{\zeta}^{\alpha \beta}(s).
\end{eqnarray*}
Now we observe that the family $f_{\zeta}^{\alpha \beta}(s)$ is
$p$-orthogonal with $q$ indices for any $(\alpha,\beta,s)$ as a
simple consequence of the $p$-orthogonality of $f$. Since $q < d$,
we can apply the induction hypothesis recalled in Remark
\ref{Remark-Induction} to obtain $$\Big\| \sum_{\zeta \in [n]^q}
f_{\zeta}^{\alpha \beta} (s) \Big\|_{L_p(\tau; S_p)} \le
\mathrm{k}_d \, p^{\frac{q(q+1)}{2}} \max_{(\varepsilon, \delta)
\in \mathbb{P}_{q}(2)} \left\{ \Big\| \sum_{\zeta \in [n]^{q}}
f_{\zeta}^{\alpha \beta}(s) \otimes e_{\pi_{\varepsilon}
(\zeta),\pi_{\delta}(\zeta)} \Big\|_p \right\}.$$ Putting it all
together, the assertion follows by Remark \ref{Remark-Converse}.
\fin

\subsection{Concluding estimates}

Now we are ready to prove Theorem \ref{Theorem-p-Orthogonality}.
First we recall that the $p$-orthogonality of $f$ can be combined
with Lemma \ref{Lemma-Injective} to drop those terms for which the
indices admit an injective projection. In other words,
\begin{eqnarray*}
\Big\| \sum_{\gamma \in \Gamma} f_{\gamma} \Big\|_{L_p(\tau)}^p &
= & \sum_{\gamma_1, \ldots, \gamma_p \in \Gamma} \tau \big(
f_{\gamma_1}^* f_{\gamma_2}^{} f_{\gamma_3}^* f_{\gamma_4}^{}
\cdots f_{\gamma_{p-1}}^* f_{\gamma_p}^{} \big) \\ & = & (-1)^d
\sum_{\sigma_1 > \dot{0}} \cdots \sum_{\sigma_d > \dot{0}} \Big[
\prod_{k=1}^d \mu(\dot{0},\sigma_k) \Big] \Psi(\sigma_1, \ldots,
\sigma_d).
\end{eqnarray*}
On the other hand, let us write
\begin{eqnarray*}
\mathsf{A} & = & \Big\| \sum_{\gamma \in \Gamma} f_{\gamma}
\Big\|_{L_p(\tau)}, \\ \mathsf{B} & = & \Big\| \sum_{i_1, \ldots,
i_d=1}^n \lambda(g_{i_1}) \otimes \cdots \otimes \lambda(g_{i_d})
\otimes f_{i_1 \cdots i_d} \Big\|_{L_p(\tau_d \otimes \tau)}, \\
\mathsf{C} & = & \max_{(\alpha, \beta) \in \mathbb{P}_d(2)}
\left\{ \Big\| \sum_{\gamma \in \Gamma} f_{\gamma} \otimes
e_{\pi_{\alpha}(\gamma),\pi_{\beta}(\gamma)} \Big\|_{L_p(\tau;
S_p)} \right\}.
\end{eqnarray*}
Let us write $\delta = (\sigma_1, \sigma_2, \ldots, \sigma_d)$ and
let $r(\delta)$ be the number of common singletons. That is,
$r(\delta)$ coincides with the cardinality of $\mathrm{B}_d$.
Then, Lemma \ref{Lemma-Factorization} and H\"{o}lder's inequality
provide the following estimate
\begin{eqnarray*}
\big| \Psi(\sigma_1, \ldots \sigma_d) \big| & \le & \prod_{s=1}^p
\|\mathsf{F}_s\|_p \\ & \le & \mathsf{A}^{r(\delta)}
\prod_{q=0}^{d-1} \big( \mathrm{k}_d p^{\frac{q(q+1)}{2}}
\mathsf{B} \big)^{|\mathrm{B}_q|} \\ & \le & \Big[
\prod_{q=0}^{d-1} p^{\frac{q(q+1)}{2} |\mathrm{B}_q|} \Big]
\mathsf{A}^{r(\delta)} \big( \mathrm{k}_d \mathsf{C}
\big)^{p-r(\delta)}.
\end{eqnarray*}
Notice that $\mathsf{B} \le \mathrm{k}_d \, \mathsf{C}$ by
inequality (\ref{Equation-Iteration}). Now, recalling that
$$\sum_{q=0}^{d-1} \frac{q(q+1)}{2} |\mathrm{B}_q| \le
\frac{d(d-1)}{2} \sum_{q=0}^{d-1} |\mathrm{B}_q| =
\frac{d(d-1)}{2} (p - r(\delta)),$$ we obtain the following
estimate for $\Psi$ $$\big| \Psi(\sigma_1, \ldots \sigma_d) \big|
\le \mathsf{A}^{r(\delta)} \Big( \mathrm{k}_d \,
p^{\frac{d(d-1)}{2}} \mathsf{C} \Big)^{p-r(\delta)}.$$ Putting it
all     together, we get $$\mathsf{A}^p \le \sum_{\sigma_1
> \dot{0}} \cdots \sum_{\sigma_d > \dot{0}} \Big[ \prod_{k=1}^d
|\mu(\dot{0},\sigma_k)| \Big] \mathsf{A}^{r(\delta)} \Big(
\mathrm{k}_d \, p^{\frac{d(d-1)}{2}} \mathsf{C}
\Big)^{p-r(\delta)}.$$ Since $\sigma_k > \dot{0}$ for all $k$, we
have $0 \le r(\delta) \le p-1$. Therefore, we can write
$$\mathsf{A}^p \le \sum_{r=0}^{p-1} \varphi_r \, \mathsf{A}^r
\Big( \mathrm{k}_d \, p^{\frac{d(d-1)}{2}} \mathsf{C}
\Big)^{p-r},$$ with $\varphi_r$ given by $$\varphi_r =
\sum_{\delta_0: \, r(\delta_0) = r} \, \prod_{k=1}^d
|\mu(\dot{0},\sigma_k)|.$$ The zero subindex in $\delta$ is chosen
to denote that the sum is taken over the set of $\delta_0 =
(\sigma_1, \sigma_2, \ldots, \sigma_d)$ such that $\sigma_k >
\dot{0}$ for all $k$. Ignoring that restriction and applying Lemma
\ref{Lemma-Mobius}, we easily get $$\varphi_r \le \sum_{\delta: \,
r(\delta) \ge r} \, \prod_{k=1}^d |\mu(\dot{0},\sigma_k)| =
{{p}\choose{r}} \prod_{k=1}^d \sum_{\sigma_k \in \mathbb{P}_{p-r}}
|\mu(\dot{0},\sigma_k)| = {{p}\choose{r}} (p-r)!^d.$$ In
particular, we obtain
\begin{eqnarray} \label{Equation-Binomial}
\mathsf{A}^p & \le & \sum_{r=0}^{p-1} {{p}\choose{r}} (p-r)!^d
\mathsf{A}^r \Big( \mathrm{k}_d \, p^{\frac{d(d-1)}{2}} \mathsf{C}
\Big)^{p-r} \\ \nonumber & \le & \sum_{r=0}^{p-1} {{p}\choose{r}}
(p-r)! \mathsf{A}^r \mathsf{D}^{p-r},
\end{eqnarray}
where $\mathsf{D}$ has the form $$\mathsf{D} = \mathrm{k}_d \Big[
\sup_{0 \le r \le p-2} (p-r)!^{\frac{d-1}{p-r}} \Big]
p^{\frac{d(d-1)}{2}} \mathsf{C} \le \mathrm{k}_d
p^{\frac{(d+2)(d-1)}{2}} \mathsf{C}.$$ The last inequality follows
easily from Stirling's formula. Now, we conclude by applying the
same arguments as in \cite{P1}. More concretely, proceeding as in
Sublemma 2.3 of \cite{P1}, we obtain
\begin{equation} \label{Equation-Estimation-AD}
\mathsf{A} \le 2 p \, \mathsf{D} \le \mathrm{k}_d \,
p^{\frac{d(d+1)}{2}} \mathsf{C}.
\end{equation}
This estimation completes the proof of Theorem
\ref{Theorem-p-Orthogonality}. Although the proof of
(\ref{Equation-Estimation-AD}) follows from
(\ref{Equation-Binomial}) and Sublemma 2.3 of \cite{P1}, we
include the proof for completeness. If $\mathsf{A} \le p \,
\mathsf{D}$ there is nothing to prove. Hence, assume that
$\mathsf{A} > p \, \mathsf{D}$. Let us divide at both sides of
(\ref{Equation-Binomial}) by $\mathsf{A}^p$ and let us write $z =
\mathsf{D} / \mathsf{A}$, so that $p z < 1$ and $$1 \le
\sum_{r=0}^{p-1} {{p}\choose{r}} (p-r)! z^{p-r}.$$ Then, we have
\begin{eqnarray*}
2 & \le & \sum_{r=0}^{p-1} {{p}\choose{r}} \int_0^{\infty}
(zx)^{p-r} e^{-x} \, dx + 1 \\ & = & \int_0^{\infty} \big[ (1 + z
x)^p - 1 \big] e^{-x} \, dx + 1 \\ & = & \int_0^{\infty} (1 + z
x)^p e^{-x} \, dx \\ & \le & \int_0^{\infty} \exp(p z x - x) \,
dx.
\end{eqnarray*}
Since $p z < 1$, we conclude $2 \le (1 - p z)^{-1}$ and $z^{-1}
\le 2p$ as desired.

\begin{remark}
\textnormal{Let us look for a moment what happens with Theorem
\ref{Theorem-p-Orthogonality} when the von Neumann algebra
$\mathcal{M}$ is commutative, so that we can think of $L_p(\tau)$
as $L_p(\mu)$ for some probability measure $\mu$. As was recalled
in \cite{P1}, if we are given a $p$-orthogonal family $f_1, f_2,
\ldots, f_n$ with one index in $L_p(\mu)$, then we obtain the
natural analog of Burkholder-Gundy inequality for a martingale
difference sequence. Namely, we have $$\Big( \int_{\Omega} \Big|
\sum_{k=1}^n f_k(\omega) \Big|^p d\mu(\omega) \Big)^{1/p} \le
\frac{3 \pi}{2} p \, \Big( \int_{\Omega} \Big[ \sum_{k=1}^n
|f_k(\omega)|^2 \Big]^{p/2} d\mu(\omega) \Big)^{1/p}.$$ In the
general case, Theorem \ref{Theorem-p-Orthogonality} provides the
following inequality $$\Big( \int_{\Omega} \Big| \sum_{\gamma \in
\Gamma} f_{\gamma}(\omega) \Big|^p d \mu(\omega) \Big)^{1/p} \le
\mathrm{k}_d p^{\frac{d(d+1)}{2}} \Big( \int_{\Omega} \Big[
\sum_{\gamma \in \Gamma} |f_{\gamma}(\omega)|^2 \Big]^{p/2}
d\mu(\omega) \Big)^{1/p}.$$}
\end{remark}

\section{Two examples}
\label{Section-Examples}

We conclude this paper with two examples of multi-indexed
$p$-orthogonal sums. The first one came out during the preparation
of \cite{PP} and was the motivation of this work. It provides a
generalization of the notion of $p$-dissociate subset of a
discrete group. The second provides an analog of the
non-commutative Burkholder-Gundy inequalities for multi-indexed
martingale difference sequences.

\subsection{Multi-indexed $p$-dissociate sets}

Let $\mathrm{G}$ be a discrete group with identity element $e$. A
subset $\Lambda = \big\{t_1, t_2, \ldots, t_n \big\}$ of
$\mathrm{G}$ is called a \emph{$p$-dissociate set} if for any
injective function $h: \{1,2, \ldots,p\} \rightarrow \{1,2,
\ldots, n\}$, the following non-cancellation property holds
$$t_{h(1)}^{-1} t_{h(2)}^{} t_{h(3)}^{-1} t_{h(4)}^{} \cdots
t_{h(p-1)}^{-1} t_{h(p)}^{} \neq e.$$ In a similar way, let
$\Gamma$ be as above and let $\Lambda = \big\{t_{\gamma}: \,
\gamma \in \Gamma \big\}$ be a subset of $\mathrm{G}$ indexed by
$\Gamma$. Then, we shall say that $\Lambda$ is a
\emph{$p$-dissociate set with $d$ indices} if the same
non-cancellation property is satisfied whenever the function $h:
\{1,2, \ldots,p\} \rightarrow \Gamma$ has an injective projection.
Let $\lambda_{\mathrm{G}}$ be the left regular representation of
$\mathrm{G}$ and let us denote by $\tau_{\mathrm{G}}$ the natural
trace on $\lambda_{\mathrm{G}}(\mathrm{G})''$. Then, since
$\tau_{\mathrm{G}}(\lambda(t))$ vanishes unless the element $t$ is
the identity $e$, it is clear that $$f = \Big\{ a_{\gamma} \otimes
\lambda_{\mathrm{G}}(t_{\gamma}) \, \Big| \ \gamma \in \Gamma
\Big\}$$ is $p$-orthogonal in $L_p(\tau \otimes
\tau_{\mathrm{G}})$ with $d$ indices for any given function $a:
\Gamma \rightarrow L_p(\tau)$. As was pointed out in \cite{P1},
these kind of sets can be used to obtain the classical
Littlewood-Paley inequalities from the case of $1$-indexed
families in Theorem \ref{Theorem-p-Orthogonality}. A remarkable
$p$-dissociate set with $d$ indices is provided by the free group
$\Free_n$ with $n$ generators $g_1, g_2, \ldots, g_n$. Indeed, let
us consider the set $$\Lambda = \Big\{ g_{i_1} g_{i_2} \cdots
g_{i_d} \, \Big| \ 1 \le i_k \le n \Big\}.$$ By the freeness of
the generators, it is not difficult to check that $\Lambda$ is
$p$-dissociate with $d$ indices for any even exponent $p$. In
particular, given any collection of operators $\mathcal{A} =
(a_{\gamma})_{\gamma \in \Gamma}$ indexed by $\Gamma$, the family
$$\mbox{\emph{\textbf{f}}}_{i_1i_2 \cdots i_d} = a_{i_1i_2 \cdots
i_d} \otimes \lambda(g_{i_1}g_{i_2} \cdots g_{i_d})$$ is
$p$-orthogonal in $L_p(\tau \otimes \tau_{\mathrm{G}})$ with $d$
indices and Theorem \ref{Theorem-p-Orthogonality} gives $$\Big\|
\sum_{\gamma \in \Gamma} \mbox{\emph{\textbf{f}}}_{\gamma}
\Big\|_{L_p(\tau \otimes \tau_{\mathrm{G}})} \le \mathrm{k}_d \,
p^{\frac{d(d+1)}{2}} \max_{(\alpha, \beta) \in \mathbb{P}_d(2)}
\left\{ \Big\| \sum_{\gamma \in \Gamma}
\mbox{\emph{\textbf{f}}}_{\gamma} \otimes e_{\pi_{\alpha}(\gamma),
\pi_{\beta}(\gamma)} \Big\|_{L_p(\tau \otimes
\tau_{\mathrm{G}};S_p)} \right\}.$$ However, the following
equivalence of norms $$\max_{(\alpha, \beta) \in \mathbb{P}_d(2)}
\left\{ \Big\| \sum_{\gamma \in \Gamma}
\mbox{\emph{\textbf{f}}}_{\gamma} \otimes e_{\pi_{\alpha}(\gamma),
\pi_{\beta}(\gamma)} \Big\|_p \right\} \simeq \max_{0 \le k \le d}
\left\{ \Big\| \sum_{\mathrm{I} \in [n]^k} \sum_{\mathrm{J} \in
[n]^{d-k}} a_{\mathrm{IJ}}^{} \otimes e_{\mathrm{I},
\mathrm{J}}^{} \Big\|_p \right\},$$ holds for any exponent $2 \le
p \le \infty$ and with constants depending only on $d$. The reader
is referred to \cite{PP} for the proof of this fact. Moreover, the
main result in \cite{PP} claims that the inequality above holds
with constants independent on $p$ and the same happens for the
reverse inequality. In summary, we have $$\Big\| \sum_{\gamma \in
\Gamma} \mbox{\emph{\textbf{f}}}_{\gamma} \Big\|_{L_p(\tau \otimes
\tau_{\mathrm{G}})} \simeq \max_{0 \le k \le d} \left\{ \Big\|
\sum_{\mathrm{I} \in [n]^k} \sum_{\mathrm{J} \in [n]^{d-k}}
a_{\mathrm{IJ}}^{} \otimes e_{\mathrm{I}, \mathrm{J}}^{}
\Big\|_{L_p(\tau; S_p)} \right\},$$ with constants depending only
on $d$. This is the second example we meet in this paper for which
the inequality in Theorem \ref{Theorem-p-Orthogonality} for
multi-indexed $p$-orthogonal sums turns out to be an equivalence
of norms, with constants depending only on $d$. The first was
given in Section \ref{Section-Iteration}, see
(\ref{Equation-Iteration}) and Remark \ref{Remark-Converse}. In
the next paragraph, we analyze one more example of this kind.

\subsection{Multi-indexed martingale difference sequences}

Let us consider a von Neumann algebra $\mathcal{M}$ with a
faithful, normal trace $\tau$ satisfying $\tau(1)=1$. For each $1
\le k \le d$, let us consider a filtration $\mathcal{M}_1(k),
\mathcal{M}_2(k), \ldots, \mathcal{M}_n(k)$ of $\mathcal{M}$. A
family $\mbox{\emph{\textbf{f}}} =
(\mbox{\emph{\textbf{f}}}_{\gamma})_{\gamma \in \Gamma}$ of random
variables in $L_p(\tau)$ will be called a \emph{martingale
difference sequence with $d$ indices} if the following condition
holds for all $k = 1,2, \ldots, d$
$$\mbox{\emph{\textbf{f}}}_{\gamma} =
\mathbb{E}_{\mathcal{M}_{i_k}(k)} \left( h_{\gamma_k}^k \right) -
\mathbb{E}_{\mathcal{M}_{i_k-1}(k)} \left( h_{\gamma_k}^k
\right).$$ Here, $\gamma_k = (i_1, \ldots, \widehat{i}_k, \ldots,
i_d)$ for $\gamma = (i_1, \ldots, i_d)$ with $\widehat{i}_k$
meaning deletion of $i_k$ and each $h^k$ is a $(d-1)$-indexed
family in $L_p(\tau)$. In other words, we require
$\mbox{\emph{\textbf{f}}}$ to be a martingale difference sequence
when looking at each component $\pi_k(\gamma) = i_k$ of the index
set $\Gamma$. Notice that we allow different filtrations for each
component. Such a construction again leads to a $p$-orthogonal
family with $d$ indices. Namely, let us assume that the $k$-th
projection of $h: \{1,2, \ldots, p\} \rightarrow \Gamma$ is
injective. Then we consider the largest value $m_k$ of $\pi_k
\circ h$ and we take conditional expectation of index $m_k - 1$
with respect to the filtration $\mathcal{M}_1(k),
\mathcal{M}_2(k), \ldots, \mathcal{M}_n(k)$ so that $$\tau \left(
\mbox{\emph{\textbf{f}}}_{h(1)}^* \mbox{\emph{\textbf{f}}}_{h(2)}
\cdots \mbox{\emph{\textbf{f}}}_{h(p-1)}^*
\mbox{\emph{\textbf{f}}}_{h(p)} \right) = \tau \left(
\mathbb{E}_{\mathcal{M}_{m_k-1}(k)} \left[
\mbox{\emph{\textbf{f}}}_{h(1)}^* \mbox{\emph{\textbf{f}}}_{h(2)}
\cdots \mbox{\emph{\textbf{f}}}_{h(p-1)}^*
\mbox{\emph{\textbf{f}}}_{h(p)} \right] \right) = 0.$$ Here
Theorem \ref{Theorem-p-Orthogonality} also admits a converse so
that we get an equivalence $$\Big\| \sum_{\gamma \in \Gamma}
\mbox{\emph{\textbf{f}}}_{\gamma} \Big\|_{L_p(\tau)} \simeq
\max_{(\alpha,\beta) \in \mathbb{P}_d(2)} \left\{ \Big\| \sum_
{\gamma \in \Gamma} \mbox{\emph{\textbf{f}}}_{\gamma} \otimes
e_{\pi_{\alpha}(\gamma), \pi_{\beta}(\gamma)} \Big\|_{L_p(\tau;
S_p)} \right\}.$$ However, in contrast with the previous
paragraph, in this case the constants depend on $d$ and $p$. This
equivalence can be regarded as the version of Burkholder-Gundy
inequalities for multi-indexed martingale difference sequences. In
fact, it can be proved without the aid of Theorem
\ref{Theorem-p-Orthogonality}. Namely, it follows easily by
iterating the (non-commutative) Khintchine inequality (with
Rademacher functions instead of free generators) and applying
repeatedly the UMD property of $L_p(\tau)$, which follows itself
by the (non-commutative) Burkholder-Gundy inequalities. The reader
is referred to the papers \cite{P1,PX,R} for more on this.

\

\noindent \textsc{Acknowledgment}. I wish to thank G. Pisier for
introducing me to this subject. I also want to thank the anonymous
referee for a very careful reading of this paper.

\bibliographystyle{amsplain}

\end{document}